\documentclass[11pt,reqno]{amsart}

\usepackage{a4wide}
\usepackage{hyperref, amssymb, amsmath, enumerate, verbatim, color, accents, bbm}
\usepackage[toc,page]{appendix}
\usepackage[dvipsnames]{xcolor}
\usepackage{amsmath, amsthm, amssymb, stmaryrd}

\hypersetup{pdfstartview={FitH}}

\newtheorem{thm}{Theorem}[section]
\newtheorem{prop}[thm]{Proposition}
\newtheorem{lem}[thm]{Lemma}
\newtheorem{cor}[thm]{Corollary}

\theoremstyle{definition}
\newtheorem{definition}[thm]{Definition}

\theoremstyle{remark}

\numberwithin{equation}{section}

\newcommand{\RR}{\mathbb{R}}

\newcommand{\ft}{\widetilde{f}}
\newcommand{\gt}{\widetilde{g}}

\newcommand{\func}[3][f]{#1\colon #2\rightarrow #3}

\newcommand{\ZZ}{\mathbb{Z}}
\newcommand{\E}{\mathbb{E}}

\let\OLDthebibliography\thebibliography
\renewcommand\thebibliography[1]{
  \OLDthebibliography{#1}
  \setlength{\parskip}{0pt}
  \setlength{\itemsep}{0pt plus 0.3ex}
}

\usepackage{mathtools}

\usepackage{multirow, array}
\usepackage{placeins}

\usepackage{caption} 
\newcommand*\wrapletters[1]{\wr@pletters#1\@nil}
\def\wr@pletters#1#2\@nil{#1\allowbreak\if&#2&\else\wr@pletters#2\@nil\fi}

\def \cP {\mathcal P}

\newenvironment{proofof}[1]{\indent{\scshape Proof of #1}:~~}{\qed}
\newenvironment{proofcap}{\indent{\scshape Proof}:~~}{\qed}

\begin{document}

\title{Szemer\' edi's theorem in the primes}

\author{Luka Rimani\'{c} and Julia Wolf}

\date{\today}

\maketitle

\begin{center}\small\emph{In memory of Kevin Henriot.}\end{center}

\begin{abstract}
Green and Tao famously proved in 2005 that any subset of the primes of fixed positive density contains arbitrarily long arithmetic progressions. Green had previously shown that in fact any subset of the primes of relative density tending to zero sufficiently slowly contains a 3-term term progression. This was followed by work of Helfgott and de Roton, and Naslund, who improved the bounds on the relative density in the case of 3-term progressions. The aim of this note is to present an analogous result for longer progressions by combining a quantified version of the relative Szemer\'edi theorem given by Conlon, Fox and Zhao with Henriot's estimates of the enveloping sieve weights.
\end{abstract}

\section{Introduction}
Let $r_k (N)$ denote the maximal size of a subset of $[N]:=\{1,2,\ldots,N\}$ not containing any non-trivial $k$-term arithmetic progressions, and let $r_{k} (\cP_N)$ denote the maximal size of a subset of the set $\mathcal{P}_N:=\mathcal{P}\cap[N]$ of primes less than $N$ not containing any non-trivial $k$-term arithmetic progressions. Define the corresponding critical densities by 
\[
\alpha_k (N) := \frac{r_k(N)}{N} \hspace{3mm}\text{and}\hspace{3mm} \alpha_{k} (\cP_N) := \frac{r_k (\cP_N)}{|\cP_N|}, 
\]
respectively. The current state of the art in the integers is
\begin{align*}
\alpha_3 (N) &\ll (\log N)^{-1+o(1)} \hspace{5mm}\text{\cite{Bloom}},\\
\alpha_4 (N) &\ll (\log N)^{-c} \hspace{3mm} \text{for some }c>0\hspace{3mm}\text{\cite{GreenTao4}},\\
\alpha_k (N) &\ll (\log \log N)^{-2^{-2^{k+9}}} \text{for all } k\geq 5 \hspace{3mm}\text{\cite{Gowers}}.
\end{align*}
Regarding relative density in the primes, the current record for progressions of length 3 is
\[
\alpha_3 (\cP_N) \ll (\log \log N)^{-1+o(1)},
\]
arrived at through a series of articles by Green \cite{Green}, Helfgott and de Roton \cite{HelfgottdeRoton}, and finally Naslund \cite{Naslund}. Henriot \cite{Henriot} extended Naslund's result to all linear systems of complexity one, and our result relies crucially on the optimised estimates of the enveloping sieve weights he gave in this paper. The aim of this note is to extend these results to longer arithmetic progressions as follows.

\begin{thm}\label{main_thm}
For every $k\geq 4$ there exists $c_k >0$ such that 
\begin{equation}\label{main_eqn}
\alpha_k (\cP_N) \ll_k \alpha_k \left( (\log \log N)^{c_k}\right).
\end{equation}
\end{thm}

To put this in perspective, Theorem \ref{main_thm} and the aforementioned results in the integers yield
\begin{align*}
\alpha_4 (\cP_N) &\ll (\log \log \log N)^{-c} \text{ for }c>0 \text{ determined by \cite{GreenTao4}},\\
\alpha_k (\cP_N) &\ll (\log \log \log \log N)^{-2^{-2^{k+9}}} \text{for all }k\geq 5.
\end{align*}
For comparison, in their celebrated work on long arithmetic progressions in the primes Green and Tao \cite{GreenTao} obtained bounds on $\alpha_{k} (\cP_N)$ of the form $(\log_{(7)} N)^{-c_k}$ whenever $k\geq 4$, where $\log_{(s)}N$ denotes the $s$-fold iterated logarithm of $N$.

The proof of Theorem \ref{main_thm} naturally splits into two parts, each of which contributes one logarithm to the bound in (\ref{main_eqn}). The first part consists of a relative Szemer\'edi theorem as proved by Conlon, Fox and Zhao \cite{ConlonFoxZhao}, which we quantify for our purposes in Section \ref{sectionrel} (some details have been relegated to the appendix). In the second part we make use of Henriot's \cite{Henriot} optimised estimates of the usual sieve weights associated with the primes. 

\section{A Quantitative relative Szemer\'edi theorem}\label{sectionrel}

At the heart of the aforementioned relative Szemer\'edi theorem lies the idea, already present in \cite{GreenTao} and elaborated on in \cite{TaoZiegler}, \cite{GowersDec}, \cite{Reingoldetal} and finally \cite{CFZRelative}, that one can deal with an unbounded function in a pseudorandom setting by approximating it with a bounded function while preserving its density and the count of arithmetic progressions. 

Let us first recall Szemer\'edi's theorem in the dense setting, which can--using Varnavides' averaging trick--be rephrased in the following weighted form.
\begin{thm}[Szemer\'edi's theorem, dense setting]\label{sz}
Suppose that $k \geq 3$ is an integer and let $\alpha >0$. Then every $\func{\ZZ_N}{[0,1]}$ with $\E f \geq \alpha$ satisfies
\begin{equation}
\E_{x,d \in \ZZ_N} [f(x)f(x+d) \ldots f(x+(k-1)d) ] \gg_k \left( \alpha_{k}^{-1} (\alpha /2) \right)^{-2},
\end{equation}
where the function $\alpha_k:\mathbb{N}\rightarrow [0,1]$ is defined as in the introduction, and $\alpha_{k}^{-1}$ denotes its inverse. 
\end{thm}

In order to extend this result to a sparse setting, we shall need to impose a pseudorandomness condition on the function in question. We shall use the following arithmetic version of a definition given by Conlon, Fox and Zhao \cite[Definition 2.2]{CFZRelative}, 
 which has the additional feature of measuring the speed of convergence. 

\begin{definition}[$(k,\delta)$-linear forms condition, arithmetic setting]\label{def_lfc_arithmetic}
Let $k\geq 2$ and $\delta >0$. We say that a function $\func[\nu]{\ZZ_N}{[0,\infty)}$ satisfies the $(k,\delta)$-\textit{linear forms condition} (or $(k,\delta)$-\textit{LFC} in short) if
\[
\left| \E_{x_{1}^{(0)},x_{1}^{(1)},\ldots, x_{k}^{(0)}, x_{k}^{(1)} \in \ZZ_N} \left[ \prod_{j=1}^{k} \prod_{\omega \in \{0,1\}^{[k]\setminus \{j\}}} \nu \left( \sum_{i=1}^{k} (j-i) x_{i}^{(\omega_i)} \right)^{n_{j,\omega}} \right]  - 1 \right| \leq \delta,
\]
for any choice of exponents $n_{j,\omega} \in \{ 0,1\}$.
\end{definition}

We are now ready to state the main result of this section.
\begin{thm}[Quantitative relative Szemer\'edi theorem]\label{relativesz}
Suppose that $k \geq 3$ is an integer and that $\func[\nu]{\ZZ_N}{[0,\infty)}$ satisfies the $(k,\delta)$-linear forms condition. Then for all $\alpha >0$ there exists $c_{k}'>0$ such that the following holds. 

If $\func{\ZZ_N}{[0,\infty)}$ is a function such that $0\leq f\leq \nu$ and $\E f \geq \alpha$, then
\begin{equation}\label{szeq}
\E_{x,d \in \ZZ_N} \left[ f(x)f(x+d) \ldots f(x+(k-1)d) \right] \gg_k \frac{1}{\left( \alpha_{k}^{-1} (\alpha /2) \right)^{2}} - \frac{1}{\log^{c_{k}'} (1/\delta)}.
\end{equation}
\end{thm}

A qualitative version of Theorem \ref{relativesz} was given by Conlon, Fox and Zhao \cite[Theorem 2.4]{CFZRelative}
, so the only novelty here is the explicit error term with respect to the speed of convergence of the linear forms condition. 

In the remainder of this section we summarise the main steps in the proof of Theorem \ref{relativesz}. To begin with we need one further definition.

\begin{definition}For any positive integer $r$ and any function $\func{\ZZ_N}{\RR}$ we define the \textit{cut norm} of $f$ by
\begin{equation}\label{cutnormf}
\|f\|_{\Box,r} := \sup | \E_{x_1, \ldots, x_{r} \in \ZZ_N }  f(x_1+ \ldots+ x_{r}) \prod_{j \in [r]}1_{A_j}(x_{-j}) |,
\end{equation}
where the supremum is taken over all $A_1,\ldots,A_r \subseteq \ZZ_{N}^{r-1}$, and $x_{-j}$ stands for the $(r-1)$-tuple $(x_1,\ldots, x_{j-1}, x_{j+1}, \ldots, x_r)\in \ZZ_{N}^{r-1}$. 
\end{definition}
In Appendix \ref{app_control_cut} we show that the $(k,\delta)$-LFC implies control in the cut norm, a qualitative version of which is implicit in \cite{CFZRelative,ConlonFoxZhao}.

\begin{cor}\label{nu_uniform}
Suppose that $\nu$ satisfies the $(k,\delta)$-linear forms condition. Then
\[
\| \nu - 1 \|_{\Box}  \leq 2 \delta^{1/2^{k-1}}.
\]
\end{cor}
This allows us to use the following \emph{dense model theorem} taken almost verbatim from \cite[Theorem 5.1]{ConlonFoxZhao}, which states that an unbounded function with a pseudorandom majorant can be well approximated by a 1-bounded function of the same density such that their difference behaves well with respect to inner products with certain test functions.

\begin{thm}[Dense model theorem]\label{densemodel}
There exists an absolute constant $C>0$ such that the following holds.

Let $\varepsilon >0$, and suppose that $\func[\nu]{\ZZ_N}{[0,\infty)}$ satisfies $\| \nu - 1 \|_{\Box, k-1} \leq \varepsilon$. Then for every $\func{\ZZ_N}{[0,\infty)}$ such that $0 \leq f\leq \nu$, there exists a function $\func[\ft]{\ZZ_N}{[0,1]}$ such that \[\|f-\ft \|_{\Box, k-1} \leq \log^{-1/C}(1/\varepsilon).\]
\end{thm}

The second ingredient in the proof of Theorem \ref{relativesz} is the so-called \textit{counting lemma}, which states that the count of arithmetic progressions is preserved under the conditions implied by the dense model theorem.

\begin{prop}[$(k,\delta)$-relative counting lemma, arithmetic setting]\label{counting}
Suppose that $\nu$ satisfies the $(k,\delta)$-linear forms condition and let $f,\ft$ be functions on $\ZZ_N$ such that $0\leq f \leq \nu$ and $0\leq \ft \leq 1$. If $\|f- \ft \|_{\Box, k-1} \leq \varepsilon$, then
\begin{align*}\label{countingerror}
&\left| \E_{x,d}f(x)f(x+d)\ldots f(x+(k-1)d) - \E_{x,d} \ft(x) \ft(x+d)\ldots \ft (x+(k-1)d) \right| \\
&\hspace{50mm}\ll_k \delta^{1/2^{2^k+k-2}} + \varepsilon^{1/2^{2^k-1}}.
\end{align*}
\end{prop}

A qualitative version of Proposition \ref{counting} was given in \cite[Theorem 2.17]{CFZRelative} (see also \cite[Lemma 4.1]{Zhao}), so again the novelty here is the explicit error term with respect to the speed of convergence of the $(k,\delta)$-LFC and the quality of the approximation of $f$ by $\ft$. We provide proofs of Corollary \ref{nu_uniform} and Proposition \ref{counting} in Appendix \ref{app_control_cut} and Appendix \ref{app_counting}, respectively, as they are relatively minor modifications of the proofs in \cite{CFZRelative,ConlonFoxZhao}.

\begin{proofof}{Theorem \ref{relativesz} assuming Corollary \ref{nu_uniform} and Proposition \ref{counting}}
Using Corollary \ref{nu_uniform} in conjunction with Theorem  \ref{densemodel} we find $\func[\ft]{\ZZ_N}{[0,1]}$ such that
\[
\|f-\ft \|_{\Box,k-1} \ll_k \frac{1}{\log^{1/C} (1/\delta)},
\]
where $C$ is the constant in Theorem \ref{densemodel}. Now Proposition \ref{counting} yields 
\begin{align*}\label{countingerror}
&\left| \E_{x,d}f(x)f(x+d)\ldots f(x+(k-1)d) - \E_{x,d} \ft(x) \ft(x+d)\ldots \ft (x+(k-1)d) \right|\\
&\hspace{20mm}\ll_k \delta^{1/2^{2^k+k-2}} + \log^{-1/2^{2^k-1} C} (1/\delta) \ll_{k} \log^{-c_k'} (1/\delta),
\end{align*}
and the claim made in Theorem \ref{relativesz} now follows from Theorem \ref{sz}. 
\end{proofof}

\section{Proof of Theorem \ref{main_thm}}

With Theorem \ref{relativesz} in hand, the remainder of the proof of Theorem \ref{main_thm} is fairly standard. It is well known that the primes behave random-like once one removes their bias with respect to residue classes modulo small primes. This is achieved using the so-called $W$-trick, which reduces Theorem \ref{main_thm} to the following.

\begin{thm}[Main theorem in $W$-tricked primes]\label{main_w}
There exists $c_k >0$ such that the following holds. 

Let $N'$ be a large integer and let $\omega := c_0 \log N'$ for some $c_0 \in [1/4,1/2]$. Let $W := \prod_{p\leq \omega}p$, and let $b$ be a positive integer coprime to $W$. Suppose that $B\subset [N']$ is a set free of $k$-term arithmetic progressions such that $b+W\cdot B \subset \cP$ and
\[
|B| \geq \alpha \frac{W}{\phi (W)} \frac{N'}{\log N'}.
\]
Then 
\[
\alpha \ll_k \alpha_{k} \left( ( \log \log N')^{c_k} \right).
\]
\end{thm}
\begin{proofof}{Theorem \ref{main_thm} assuming Theorem \ref{main_w}}
Let $N$ be a sufficiently large  integer and let $A \subset \cP \cap [N]$ be a subset of relative density $\alpha$ free of $k$-term progressions. We may assume that $\alpha \gg N^{-1/4}$ for otherwise the claim is trivial. Let $\omega := \frac14 \log N$. It is not difficult to see that for such a choice of $\omega$ one has $W=N^{1/4+o(1)}$. Define $N' := \lfloor N / W \rfloor = N^{3/4+o(1)}$. By an averaging argument (see for example Lemma 2.1 in \cite{HelfgottdeRoton}) there exists a positive integer $b$ coprime to $W$ such that the set $B := \{ n \in [N'] \colon b+Wn \in A\}$ satisfies
\[
|B| \gg \alpha \frac{W}{\phi (W)} \frac{N'}{\log N'}.
\]
Since $\omega \sim \frac13 \log N'$ as $N\rightarrow\infty$ and $B$ inherits the property of not containing any non-trivial $k$-term arithmetic progressions (by translation-invariance and homogeneity), we can use Theorem \ref{main_w} to conclude that
\[
\alpha \ll_k \alpha_k \left( (\log \log N')^{c_k} \right) \ll_k \alpha_k \left( (\log \log N)^{c_k} \right),
\]
which completes the proof of Theorem \ref{main_thm}.
\end{proofof}

It thus suffices to prove the main theorem in the $W$-tricked primes. In order to use Theorem \ref{relativesz} one needs to construct a majorant for the primes that satisfies the $(k,\delta)$-LFC. The existence of such a pseudorandom majorant was established already in \cite{GreenTao} using ingredients from \cite{GoldstonPintzYildirim}. It is here that we make crucial use of the simplified and optimised estimates from \cite{Henriot}, allowing us to choose $w$ as large as $c_0 \log N'$. 

\begin{proofof}{Theorem \ref{main_w}}
By Bertrand's postulate one can find a prime $M \in [2N', 4N']$. We define
\begin{equation}\label{def_lambda}
\lambda_{b,W} (n) := \frac{\phi (W)}{W} (\log N') 1_{[N']} (n) 1_{\cP}(b+Wn),
\end{equation}
to avoid possible wrap-around issues. Bringing the indicator function of $B\subseteq [N']$ in line with these weights, we define
\[
f_B (n) := \frac{M}{N'} \frac{\phi (W)}{W} (\log N') 1_B (n).
\]
Note that $\E_{[M]} f_B \geq \alpha$, and $0 \leq f_B \ll \lambda_{b,W}$ since $b+W\cdot B \subset \cP$. 

By Proposition 6.2 in \cite{Henriot}, since $M$ is large enough and equal to $N'$ up to a constant factor, there exists a pseudorandom majorant $\func[\nu]{\ZZ_M}{\RR^{+}}$ satisfying the $(k,\delta)$-LFC with $\delta \ll_k (\log N')^{-1+o(1)}$, as well as the inequality
\[
0 \leq f_{B} \ll \lambda_{b,W} \ll_k \nu.
\]
Since $B$ does not contain any non-trivial $k$-term arithmetic progressions, the left-hand side of (\ref{szeq}) with $N=M$ and $f=f_B$ is $O(1/N')$, which is easily seen to be negligible. It therefore follows from Theorem \ref{relativesz} that
\[
\frac{1}{\alpha_{k}^{-1} (\alpha/2)} \ll_k  \frac{1}{\log^{c_k'/2} (1/\delta)} \ll_k \frac{1}{(\log \log N')^{c_k}}
\]
for some constant $c_k$. This implies Theorem \ref{main_w} as $\alpha_k$ is a non-increasing function.
\end{proofof}

\section{Remarks}
We do not see any fundamental obstruction to extending the result in this paper to more general systems of translation-invariant linear equations, but shall not attempt to do so here. It would arguably be of greater interest to find a more direct approach for longer progressions along the lines of \cite{HelfgottdeRoton} and \cite{Naslund}, where the bounded function $\widetilde{f}$ is constructed explicitly.

\section{Acknowledgements}
The authors would like to thank Ben Green, Terence Tao and Aled Walker for several useful discussions. The second author is grateful to R\' egis de la Bret\`eche for an invitation to Paris in 2016 which prompted this work.

\appendix 
\section{Hypergraph notation}\label{app}

In the subsequent appendices we prove the results claimed in Section \ref{sectionrel} by quantifying the arguments in \cite{CFZRelative,ConlonFoxZhao}. These are given in the language of hypergraphs rather than linear systems of equations, and since the notation in the former setting is more compact and closer to the original source of the arguments we shall use it here. The purpose of the present section is therefore to set up the hypergraph notation used in \cite{ConlonFoxZhao}, and to indicate how it relates to that used throughout the main body of the note.

In what follows, unless otherwise mentioned, we always suppose that $k$ is a fixed positive integer and that the hypergraphs in question are weighted $(k-1)$-uniform  $k$-partite hypergraphs on $X:= X_1 \cup X_2 \cup \ldots \cup X_k$. When applied to the arithmetic setting each set $X_i$ corresponds to a copy of $\ZZ_N$.
To keep the notation as light as possible, write 
\[X_{-i} := X_1 \times \ldots \times X_{i-1} \times X_{i+1} \times \ldots \times X_k\] 
and 
\[x_{-i} := (x_1,\ldots, x_{i-1}, x_{i+1}, \ldots, x_k ) \in X_{-i}.\]
Our weighted $(k-1)$-uniform $k$-partite hypergraph then consists of a $k$-tuple of functions $g=(g_{-i})_{i=1,\dots,k}$, with each $g_{-i}: X_{-i}\rightarrow \RR$. For two such hypergraphs $g$ and $\nu$ we say that $g \leq \nu$ whenever $g_{-i} (x_{-i}) \leq \nu_{-i} (x_{-i})$ for all $1\leq i \leq k$ and all $x_{-i} \in X_{-i}$. We shall often omit the index of the weight when it is clear from the argument it takes. For example, by $\prod_{i=1}^{k} g (x_{-i})$ we mean $\prod_{i=1}^{k} g_{-i} (x_{-i})$.

Given a weighted $r-$uniform hypergraph $h$ on $X_{1}\times \ldots \times X_r$, that is, a function $h:X_{1}\times \ldots \times X_r\rightarrow \RR$, we define the \textit{cut norm} of $h$ by
\begin{equation}\label{cutnormhg}
\| h \|_{\Box,r} := \sup | \E_{x_1 \in X_1, \ldots, x_r \in X_r } h(x_1, \ldots, x_r) \prod_{j \in [r]}1_{A_j}(x_{-j}) |,
\end{equation}
where the supremum is taken over all $A_j \subseteq X_{-j}$. Finally, given a weighted $(k-1)$-uniform $k$-partite hypergraph $g$ on vertex set $X$ as above, we define the \textit{cut norm} of $g$ by
\begin{equation}\label{cutnormmax}
\|g\|_{\Box} := \max \left\{ \|g_{-1} \|_{\Box, k-1}, \ldots, \|g_{-k}\|_{\Box, k-1} \right\}.
\end{equation}
It is not difficult to see that both $\|.\|_{\Box,r}$ and $\|.\|_{\Box}$ are indeed norms. They are related to their arithmetic counterparts in Section \ref{sectionrel} by setting $h(x_1,\ldots,x_r) = f(x_1+\ldots+x_r)$. 

With this notation we have following analogue of the $(k,\delta)$-LFC introduced in Definition \ref{def_lfc_arithmetic}.
\begin{definition}[$(k,\delta)$-linear forms condition, hypergraph setting]
Let $\delta >0$.  We say that a weighted $(k-1)$-uniform $k$-partite hypergraph $\nu$ satisfies the $(k,\delta)$-\textit{linear forms condition} ($(k,\delta)$-\textit{LFC}) if
\begin{equation}\label{LFC}
\left| \E_{x_{1}^{(0)},x_{1}^{(1)} \in X_1,\ldots, x_{k}^{(0)}, x_{k}^{(1)} \in X_k}\left[ \prod_{j=1}^{k} \prod_{\omega \in \{ 0,1\}^{[k]\setminus \{j\}}} \nu \left( \left(x^{(\omega)}\right)_{-j}\right)^{n_{j,\omega}} \right] - 1 \right| \leq \delta,
\end{equation}
for any choice of exponents $n_{j,\omega} \in \{ 0,1\}$, with $x^{(\omega)} = (x_{1}^{(\omega_1)}, \ldots, x_{k}^{(\omega_k)} )$.
\end{definition}
In other words, satisfying the $(k,\delta)$-LFC in a hypergraph amounts to containing the expected count of every subgraph of the 2-blow-up of $K_{k}^{(k-1)}$. Observe also that the $(k,\delta)$-LFC in Definition \ref{def_lfc_arithmetic} is easily recovered from the above by making the substitution
\[
 \left(x^{(\omega)}\right)_{-j} \longmapsto \sum_{i=1}^{k} (j-i) x_{i}^{(\omega_i)}.
\]

\section{Controlling the cut norm}\label{app_control_cut}

We begin by proving the following auxiliary result, which will come in useful when attempting to bound the cut norm of $\nu-1$, as well as in establishing the so-called \textit{strong linear forms condition} in Appendix \ref{app_strong}.

\begin{lem}\label{lem_LFC}
Suppose that $\nu$ satisfies the $(k,\delta)$-linear forms condition. For each $j\in [k]$ and each $\omega \in \{0,1\}^{[k]\setminus \{j \}}$, let $h_{j,\omega} \in \{ 1, \nu, \nu - 1\}$. Let $K$ denote the number of pairs $(j,\omega)$ for which $h_{j,\omega} = \nu - 1$ and suppose that $K\geq 1$. Then 
\[
S \left( \left(h_{j,\omega} \right)_{j,\omega} \right) := \E_{x_{1}^{(0)}, x_{1}^{(1)} \in X_1, \ldots x_{k}^{(0)}, x_{k}^{(1)} \in X_k} \prod_{j=1}^{k} \prod_{\omega \in \{ 0,1\}^{[k]\setminus \{j\}}} h_{j,\omega} \left( \left(x^{(\omega)}\right)_{-j}\right)
\]
satisfies the inequality
\begin{equation}
\left| S \left( \left(h_{j,\omega} \right)_{j,\omega} \right)  \right| \leq 2^K \delta .\end{equation}
\end{lem}
\begin{proofcap}
If all $h_{j,\omega}$ are equal to $1$ or $\nu$, then since $\nu$ satisfies the $(k,\delta)$-LFC, we have $S ( \left(h_{j,\omega} \right)_{j,\omega}) \in [ 1-\delta, 1+ \delta ]$. In the general case, when expanding $S ( \left(h_{j,\omega} \right)_{j,\omega})$ we get to decide whether to choose $\nu$ or $-1$ for each pair $(j,\omega)$ such that $h_{j,\omega} = \nu - 1$.  The terms that make a positive contribution are exactly those in which the number of $-1$s is even, the terms in which the number of $-1$s is odd making a negative contribution instead. Thus
\begin{align*}
\left| S \left( \left(h_{j,\omega} \right)_{j,\omega} \right)  \right| \leq \left[ {K \choose 0} + {K \choose 2} + \ldots \right] (1+\delta) -\left[ {K \choose 1} + {K \choose 3} + \ldots \right] (1-\delta) = 2^K \delta.
\end{align*}
\end{proofcap}

Instead of Corollary \ref{nu_uniform} we shall actually prove the slightly stronger statement that the $(k,\delta)$-LFC yields control in an appropriate Gowers uniformity norm, which dominates the cut norm and is defined as follows. Given an $r$-uniform hypergraph $h$ on $X_1 \times \ldots \times X_r$, define the \textit{Gowers $U^r$-norm} by
\[
\|h \|_{U^r}^{2^r} := \E_{x_{1}^{(0)},x_{1}^{(1)} \in X_1, \ldots, x_{r}^{(0)}, x_{r}^{(1)} \in X_r} \prod_{j=1}^{r} \prod_{\omega \in \{ 0,1\}^{[r]}} h \left(x^{(\omega)}\right).
\]

For $k=3$, the following statement appears in qualitative form as Lemma 6.3 in \cite{ConlonFoxZhao}.

\begin{cor}[$(k,\delta)$-LFC implies uniformity]\label{LFC_unif_app}
Suppose that $\nu$ satisfies the $(k,\delta)$-linear forms condition. Then
\[
\| \nu - 1 \|_{\Box} \leq \| \nu - 1 \|_{U^{k-1}} \leq 2 \delta^{1/2^{k-1}}.
\]
\end{cor}
\begin{proofcap}
Upon recalling that $\|.\|_{\Box}$ is the maximum over all $\|\cdot \|_{\Box, k-1}$, the first inequality is a straightforward application of the Gowers-Cauchy-Schwarz inequality \cite{Gowers} formulated for hypergraphs (see, for example, \cite[Lemma 6.2]{CFZRelative}). Note that the $U^{k-1}$-norm consists of $2^{k-1}$ factors, each equal to $\nu-1$, so that the result follows from Lemma \ref{lem_LFC} with $K=2^{k-1}$.\end{proofcap}

\section{Strong linear forms condition}\label{app_strong} 

The proof of the counting lemma in \cite{CFZRelative}, which we shall follow in Appendix \ref{app_counting}, proceeds by induction on the number of majorants  that are not identically 1. This requires us to be able to replace $\nu$ by $1$ under certain assumptions, which we shall be able to do as a result of the following lemma, known as the \textit{strong linear forms condition}.

\begin{lem}[$(k,\delta)$-strong linear forms]\label{SFC_app}
Let $\delta>0$. Suppose that $\nu$ satisfies the $(k,\delta)$-linear forms condition and let $0\leq g \leq \nu$, $0\leq \gt \leq 1$. Then
\begin{equation}\label{eq:exp}  
\left| \E_{\tiny\begin{array}{c} x_{1}\in X_1,\ldots, x_{k-1} \in X_{k-1}  \\ x_{k}^{(0)}, x_{k}^{(1)} \in X_{k} \end{array}}  \left( \nu (x_{-k}) - 1 \right) \prod_{j=1}^{k-1} \prod_{\omega \in \{ 0,1\}} h_{j,\omega} \left( \left(x^{(\omega)}\right)_{-j}\right) \right| \leq 2 (1+\delta)^{1 - 1/2^{k-1}} \delta^{1/2^{k-1}},
\end{equation}
where $x^{(\omega)}:=(x_1,\dots,x_{k-1},x_k^{(\omega)})$ and each $h_{j,\omega} \in \{ g_{-j}, \gt_{-j} \}$.
\end{lem}
It is the explicit form of the error term, given here as a function of the speed of convergence in the $(k,\delta)$-LFC, which distinguishes this statement from \cite[Lemma 6.3]{CFZRelative}.
\begin{proofcap}
We may without loss of generality suppose that each $h_{j,\omega}=g_{-j}$, as having $\gt$ instead of $g$ is in fact advantageous, allowing us to replace any upper bound of $\nu$ by $1$. Denote the expectation in (\ref{eq:exp})  by $S_k\left( g\right)$. Applying the Cauchy-Schwarz inequality for the first time, isolating the $x_{k-1}$ variable, we obtain
\begin{align*}
&\hspace{30mm}\left| S_k\left( g \right) \right|^2 \leq \left| \E_{\tiny\begin{array}{c} x_{1}\in X_1,\ldots, x_{k-2} \in X_{k-2}  \\ x_{k}^{(0)}, x_{k}^{(1)} \in X_{k} \end{array}}  \prod_{\omega \in \{ 0,1\}} g \left( x^{(\omega)}_{-(k-1)}\right)\right| \times \\
& \left|  \E_{\tiny\begin{array}{c} x_{1}\in X_1,\ldots, x_{k-2} \in X_{k-2}  \\ x_{k}^{(0)}, x_{k}^{(1)} \in X_{k} \end{array}}  \prod_{\omega \in \{ 0,1\}} g \left( x^{(\omega)}_{-(k-1)}\right)\left(\E_{x_{k-1} \in X_{k-1}} \left( \nu (x_{-k}) - 1 \right)  \prod_{j=1}^{k-2} \prod_{\omega \in \{ 0,1\}} g \left( x^{(\omega)}_{-j}\right)\right)^2 \right|.
\end{align*}
Bounding instances of $g_{-(k-1)}$ above by $\nu$, the above expression is at most
\begin{align*}
&\hspace{30mm}\left| \E_{\tiny\begin{array}{c} x_{1}\in X_1,\ldots, x_{k-2} \in X_{k-2}  \\ x_{k}^{(0)}, x_{k}^{(1)} \in X_{k} \end{array}} \prod_{\omega \in \{ 0,1\}} \nu \left( x^{(\omega)}_{-(k-1)}\right) \right| \times \\
&\left| \E_{\tiny\begin{array}{c} x_{1}\in X_1,\ldots, x_{k-2} \in X_{k-2}  \\ x_{k}^{(0)}, x_{k}^{(1)} \in X_{k} \end{array}} \prod_{\omega \in \{ 0,1\}} \nu \left( x^{(\omega)}_{-(k-1)}\right)\left(\E_{x_{k-1} \in X_{k-1}} \left( \nu (x_{-k}) - 1 \right)  \prod_{j=1}^{k-2} \prod_{\omega \in \{ 0,1\}} g \left( x^{(\omega)}_{-j}\right) \right)^2\right|.
\end{align*}
Denoting the second expectation by $S_{k-1}(g)$, the $(k,\delta)$-LFC now implies that
\[\left| S_k\left( g \right) \right|^2\leq (1+\delta) \left| S_{k-1} (g) \right|.\]
In order to bound $S_{k-1}(g)$, expand the square and use the Cauchy-Schwarz inequality on $x_{k-2}$. Continuing inductively we obtain
\[
|S_k (g)|^{2^{k-1}} \leq (1+\delta)^{2^{k-2}} |S_{k-1}(g)|^{2^{k-2}} \leq \ldots \leq (1+\delta)^{2^{k-1}-1} |S_{1}(g)|,
\]
where
\[
S_1 (g) := \E_{x_{1}^{(0)},x_{1}^{(1)} \in X_1, \ldots, x_{k}^{(0)}, x_{k}^{(1)} \in X_k}\prod_{\omega \in \{0,1\}^{[k-1]}} \left( \nu (x_{-k}^{(\omega)}) -1 \right) \prod_{j=1}^{k-1}\prod_{\omega \in \{0,1\}^{[k] \setminus \{j\}}} \nu (x_{-j}^{(\omega)})
\]
and the exponent of $(1+\delta)$ arises as the sum of $2^{k-j}$ from $j=2$ to $k$.
There are $2^{k-1}$ terms of the form $\nu - 1$, so by Lemma \ref{lem_LFC} we have that
\[
|S_{k} (g) | \leq 2^{2^{k-1}} \delta, 
\]
from which the claim easily follows.
\end{proofcap}

\section{The counting lemma}\label{app_counting}

In this section we shall prove the relative counting lemma, Proposition \ref{counting}, whose hypergraph version is of the following form.

\begin{prop}[$(k,\delta)$-relative counting lemma, hypergraph setting]\label{counting_app}
Let $\varepsilon>0$, and let $\nu$ be a weighted hypergraph satisfying the $(k,\delta)$-linear forms condition. Suppose that $0\leq g \leq \nu$, $0\leq \gt \leq 1$. If $\|g- \gt \|_{\Box} \leq \varepsilon$, then
\begin{equation}\label{countingerror}
\left| \E_{x_{1}\in X_1,\ldots, x_{k} \in X_{k}} \left[ \prod_{j=1}^{k} g(x_{-j}) -  \prod_{j=1}^{k} \gt (x_{-j}) \right] \right| \ll_k \delta^{1/2^{2^k+k-2}} + \varepsilon^{1/2^{2^k-1}}.
\end{equation}
\end{prop}

The purpose of what follows is to keep track of the dependence on $\delta$ and $\epsilon$ in the argument as given in \cite{CFZRelative}. In fact, it is not difficult to see that this dependence is polynomial, but we shall be a little more precise here (and attempt to make this paper at least somewhat self-contained).

The proof proceeds by induction on the number of $\nu_{-1}, \ldots, \nu_{-k}$ that are not identically $1$. Denote this number by $m$, and for every $m$ denote the least upper bound on the left-hand side of (\ref{countingerror}) by $C_m(\varepsilon, \delta)$. The case $m=0$ is addressed by the following statement.

\begin{prop}[Dense counting lemma, hypergraph setting]\label{denseprop_app}
Let $0\leq g, \gt \leq 1$. If $\|g- \gt \|_{\Box} \leq \varepsilon$, then
\begin{equation}\label{countingC_0}
\left| \E_{x_{1}\in X_1,\ldots, x_{k} \in X_{k}} \left[ \prod_{j=1}^{k} g(x_{-j}) -  \prod_{j=1}^{k} \gt (x_{-j}) \right] \right| \leq k\varepsilon.
\end{equation}
\end{prop}
For the sake of completeness, we include a proof of Proposition \ref{denseprop_app} for general $k$, given in \cite[Proposition 6.1]{ConlonFoxZhao} for $k=3$.

\begin{proofcap}
For $i\neq j$, let 
\[
X_{-i,j} := X_1 \times \ldots \times X_{i-1} \times X_{i+1} \times \ldots \times X_{j-1} \times X_{j+1} \times \ldots \times X_{k},
\]
and define $x_{-i,j}$ in a similar vein. Recall that $\|g\|_{\Box}$ is the maximum of $\|g_{-1}\|_{\Box, k-1}, \ldots, \|g_{-k}\|_{\Box, k-1}$, so that $\|g- \gt \|_{\Box} \leq \varepsilon$ implies that for all functions $\func[a_{j}]{X_{-j,k}}{[0,1]}$ we have
\[
\left| \E_{x_1\in X_1, \ldots, x_{k-1}\in X_{k-1}} (g(x_{-k}) - \gt (x_{-k})) a_1 (x_{-1,k}) \ldots a_{k-1}(x_{-(k-1),k})\right| \leq \varepsilon.
\]
Indeed, it is not difficult to see that this condition is equivalent to that for $\{0,1\}$ valued functions, which are given by the definition of the cut norm. 
Now for fixed $x_k \in X_k$, we can set $a_i (x_{-i,k}) := g (x_{-i})$ in the above expectation, and thus
\[
\left| \E_{x_{1}\in X_1,\ldots, x_{k} \in X_{k}}  (g(x_{-k}) - \gt (x_{-k})) \prod_{j=1}^{k-1} g(x_{-j}) \right| \leq \varepsilon.
\]
Similarly we get
\[
\left| \E_{x_{1}\in X_1,\ldots, x_{k} \in X_{k}}  \gt (x_{-k}) (g(x_{-(k-1)}) - \gt (x_{-(k-1)})) \prod_{j=1}^{k-2} g(x_{-j}) \right| \leq \varepsilon
\]
and so on. Continuing to telescope in this way, and using the triangle inequality, the claimed statement follows.
\end{proofcap}

We are now ready to prove the relative counting lemma. 

\begin{proofof}{Proposition \ref{counting_app}}
Proposition \ref{denseprop_app} implies that $C_0(\varepsilon, \delta) \leq k\varepsilon$. Suppose that we have calculated upper bounds for all $C_m (\varepsilon, \delta)$ for $m=0, 1, \ldots, M$, where $m$ denotes the number of $\nu_{-1}, \ldots, \nu_{-k}$ that are not identically $1$, as above. 
Suppose now that $m=M+1$, and without loss of generality that $\nu_{-1}$ is not identically one. We define auxiliary weighted $(k-1)$-uniform hypergraphs $\func[\nu', g', \gt']{X_{-1}}{[0,\infty)}$ by
\begin{align*}
\nu' (x_{-1}) &:= \E_{x_1 \in X_1} \left[\nu (x_{-2})\ldots \nu (x_{-k})\right], \\
g' (x_{-1}) &:= \E_{x_1 \in X_1} \left[g (x_{-2})\ldots g (x_{-k})\right], \\
\gt' (x_{-1}) &:= \E_{x_1 \in X_1} \left[\gt (x_{-2})\ldots \gt (x_{-k})\right].
\end{align*}
Unlike $\gt'$, the new functions $g'$ and $\nu'$ may not be bounded by $1$. We therefore define $g_{\wedge 1}' := \max \{ g',1\}$ and $\nu_{\wedge 1}' := \max \{ \nu',1\}$. As noted in \cite{CFZRelative}, the main idea is that $g',\nu'$ behave like dense graphs so that capping them by 1 produces only a small error. In what follows all expectations will be taken over $X_{-1} = X_2 \times \ldots \times X_k$ unless otherwise stated. First note that the $(k,\delta)$-LFC implies that
\begin{equation*}
1-\delta \leq \E \left[ \nu'\right], \E \left[ \nu'^{2}  \right] \leq 1+ \delta,
\end{equation*}
from which it follows, by Cauchy-Schwarz, that
\begin{equation}\label{var_nu}
\left( \E \left[ \left| \nu' -1 \right| \right] \right)^2 \leq \E \left[ (\nu' -1)^2 \right] \leq 3 \delta.
\end{equation}
\textbf{Claim.} We have
\begin{equation}\label{cap_small}
\| g_{\wedge 1}^{'} - \gt'\|_{\Box, k-1} \leq (3\delta)^{1/2} + C_M (\varepsilon, \delta) .
\end{equation}
\begin{proofof}{Claim}
Since $0\leq g' \leq \nu '$, we have
\begin{equation}\label{cap_nu}
0 \leq g' - g_{\wedge ,1}' = \max \{ g' - 1, 0 \} \leq \max \{ \nu' - 1, 0 \} \leq | \nu ' - 1|. 
\end{equation}
Note that for any $A_2 \subseteq X_{-1,2}, \ldots, A_k \subseteq X_{-1,k}$, we can write
\begin{align*}
&\E \left[ \left( g_{\wedge 1}' - \gt' \right) (x_{-1}) 1_{A_2} (x_3,\ldots,x_k) \ldots 1_{A_k} (x_2,\ldots,x_{k-1}) \right]\\
=&\E \left[ \left( g_{\wedge 1}' - g' \right) (x_{-1}) 1_{A_2} (x_{-1,2})\ldots 1_{A_k} (x_{-1,k}) \right] +\E \left[ \left( g' - \gt' \right) (x_{-1}) 1_{A_2} (x_{-1,2})\ldots 1_{A_k} (x_{-1,k}) \right].
\end{align*}
It follows from (\ref{cap_nu}) and (\ref{var_nu}) that the first expectation is at most $(3\delta)^{1/2}$ in magnitude. In order to estimate the second term, rewrite it as 
\begin{equation}\label{equation_claim}
\E_{x_{1},\ldots, x_{k}} \left[ 1_{A_2, \ldots, A_k} (x_{-1}) \prod_{j=2}^{k} g(x_{-j}) -  1_{A_2, \ldots, A_k} (x_{-1}) \prod_{j=2}^{k} \gt (x_{-j})\right],
\end{equation}
where \[1_{A_2, \ldots, A_k} (x_{-1}) :=  \prod_{j=2}^{k} 1_{A_j} (x_{-1,j}).\]
In this form it is easy to see that the number of factors whose majorant is not identically 1 is at most $M$, so by the inductive hypothesis (\ref{equation_claim}) is bounded by $C_M (\varepsilon, \delta)$. This completes the proof of the claim.
\end{proofof}

Returning to the expression we set out to bound, we see that
\[\E_{x_{1},\ldots, x_{k}} \left[ \prod_{j=1}^{k} g(x_{-j}) -  \prod_{j=1}^{k} \gt (x_{-j}) \right] =  \E_{x_{2},\ldots, x_{k}} \left[  g(x_{-1}) g' (x_{-1}) -  \gt (x_{-1}) \gt ' (x_{-1}) \right],\]
which in turn we can rewrite as
\[ \E \left[  g \left(g'  -  \gt' \right)  \right] +
\E \left[  \left( g - \gt  \right) \gt'  \right].\]
Recall that $0 \leq \gt \leq 1$ so the second term is at most
\begin{equation}\label{errorone}
\| g- \gt \|_{\Box} \leq \varepsilon,
\end{equation} 
as in the proof of Proposition \ref{denseprop_app}. Concerning the first term, the Cauchy-Schwarz inequality yields
\[\left( \E \left[  g\left(g'  -  \gt' \right)  \right] \right)^2 \leq \E\left[  g \left(g' -  \gt' \right)^2  \right] \E \left[  g \right] \leq  \E  \left[  \nu \left(g'  -  \gt' \right)^2  \right]  \E \left[  \nu \right],\]
which, upon expanding the square, using the $(k,\delta)$-LFC and the strong linear forms condition (Lemma \ref{SFC_app}), is bounded above by
\begin{equation}\label{errortwo}
(1+\delta) \left( 8 (1+\delta)^{1 - 1/2^{k-1}} \delta^{1/2^{k-1}} + \E \left[  \left(g'  -  \gt' \right)^2  \right]  \right).
\end{equation}
Expanding the final term further as
\[
\E \left[  ( g'  -  \gt'  )^2  \right] = \E \left[ (g' - \gt' )(g' - g_{\wedge 1}') \right] + \E \left[ (g' - \gt' )(g_{\wedge 1}' - \gt ') \right],
\]
we observe that since $0\leq g'\leq \nu'$ and $0\leq \gt'\leq 1$ the first term is bounded by
\[\E \left[ \nu' \left| \nu' -1 \right| \right] = \E \left[(\nu' -1)\left| \nu' -1 \right| \right] + \E \left[\left| \nu' -1 \right| \right] \leq 3\delta + (3\delta)^{1/2},\]
by (\ref{var_nu}). Finally, we rewrite the second term as
\begin{equation}\label{summands}
 \E \left[ (g' - \gt' )(g_{\wedge 1}' - \gt ') \right] = \E \left[ g' g_{\wedge 1}' \right] - \E \left[ g' \gt' \right] - \E \left[ \gt' g_{\wedge 1}' \right] + \E \left[ (\gt')^2\right].
\end{equation}
We claim that each of the four summands is close to $\E [(\gt ')^2]$. Indeed, we can write
\begin{equation}\label{errorthree}
\E [g' g_{\wedge 1}'] - \E [(\gt')^2] = \E_{x_1,\ldots,x_k} \left[ g_{\wedge 1}' (x_{-1}) \prod_{j=2}^{k} g(x_{-j}) - \gt' (x_{-1}) \prod_{j=2}^{k} \gt (x_{-j}) \right],
\end{equation}
and observe that the tuples $(g_{\wedge 1}',g_2,\dots,g_k)$ and $(\gt',\gt_{-2},\dots,\gt_{-k})$ satisfy the box-norm condition in Proposition \ref{counting_app} with $\varepsilon$ replaced by the upper bound in (\ref{cap_small}) while the number of majorants that are not identically 1 is at most $M$, so that by the inductive hypothesis (\ref{errorthree}) is bounded above by
\begin{equation}\label{errorfour}
C_M \left( (3\delta)^{1/2} + C_M (\varepsilon,\delta), \delta \right).
\end{equation}
Bounding the difference between $\E [(\gt ')^2]$ and each of the remaining summands in (\ref{summands}) in a similar fashion, we obtain
\begin{equation}\label{errorfive}
\E [ (g' - \gt')^2] \leq 3\delta + (3\delta)^{1/2} + 3C_M \left( (3\delta)^{1/2} + C_M (\varepsilon,\delta),\delta \right).
\end{equation}
Collecting all the error terms, namely (\ref{errorone}), (\ref{errortwo}) and (\ref{errorfive}), we see that
\begin{equation}\label{finalerror}
C_{M+1} (\varepsilon,\delta ) \ll_k \varepsilon + \left( \delta^{1/2^{k-1}} + C_M \left(3\delta^{1/2} + C_M(\varepsilon,\delta),\delta \right)\right)^{1/2}.
\end{equation}
For $M\geq 1$ this is satisfied whenever
\[
C_{M+1}(\varepsilon,\delta) \ll_k C_M (C_M(\varepsilon,\delta),\delta)^{1/2},
\]
which in turn implies the bound
\[
C_k (\varepsilon,\delta) \ll_{k} \delta^{1/2^{2^k+k-2}} + \varepsilon^{1/2^{2^k-1}}.
\]
This completes the proof of Proposition \ref{counting_app}.
\end{proofof}

\end{document}